\documentclass[12pt,leqno]{article}
\usepackage{amsfonts}
\usepackage{amsmath, amsthm, amsfonts, amssymb, color}
\usepackage{mathrsfs}
\usepackage{color}
\theoremstyle{definition}

\newcommand{\scr}[1]{\mathscr #1}
\definecolor{wco}{rgb}{0.5,0.2,0.3}

\numberwithin{equation}{section} \theoremstyle{remark}

\newcommand{\ua}{\uparrow}

\title{{\bf Exponential Contraction in Wasserstein Distances for Diffusion Semigroups with Negative Curvature}\footnote{Supported in
 part by  NNSFC (11131003, 11431014).} }
\author{
{\bf     Feng-Yu Wang  }\\
\footnotesize{   Laboratory of Mathematical and  Complex Systems,
Beijing Normal
University, Beijing 100875, China}\\
 \footnotesize{ Department of Mathematics,
Swansea University, Singleton Park, SA2 8PP, United Kingdom}\\
\footnotesize{  wangfy@bnu.edu.cn, F.-Y.Wang@swansea.ac.uk}}
\begin{document}
\allowdisplaybreaks
\def\R{\mathbb R}  \def\ff{\frac} \def\ss{\sqrt} \def\B{\mathbf
B}
\def\N{\mathbb N} \def\kk{\kappa} \def\m{{\bf m}}
\def\ee{\varepsilon}\def\ddd{D^*}
\def\dd{\delta} \def\DD{\Delta} \def\vv{\varepsilon} \def\rr{\rho}
\def\<{\langle} \def\>{\rangle} \def\GG{\Gamma} \def\gg{\gamma}
  \def\nn{\nabla} \def\pp{\partial} \def\E{\mathbb E}
\def\d{\text{\rm{d}}} \def\bb{\beta} \def\aa{\alpha} \def\D{\scr D}
  \def\si{\sigma} \def\ess{\text{\rm{ess}}}
\def\beg{\begin} \def\beq{\begin{equation}}  \def\F{\scr F}
\def\Ric{\text{\rm{Ric}}} \def\Hess{\text{\rm{Hess}}}
\def\e{\text{\rm{e}}} \def\ua{\underline a} \def\OO{\Omega}  \def\oo{\omega}
 \def\tt{\tilde} \def\Ric{\text{\rm{Ric}}}
\def\cut{\text{\rm{cut}}} \def\P{\mathbb P} \def\ifn{I_n(f^{\bigotimes n})}
\def\C{\scr C}      \def\aaa{\mathbf{r}}     \def\r{r}
\def\gap{\text{\rm{gap}}} \def\prr{\pi_{{\bf m},\varrho}}  \def\r{\mathbf r}
\def\Z{\mathbb Z} \def\vrr{\varrho} \def\ll{\lambda}
\def\L{\scr L}\def\Tt{\tt} \def\TT{\tt}\def\II{\mathbb I}
\def\i{{\rm in}}\def\Sect{{\rm Sect}}  \def\H{\mathbb H}
\def\M{\scr M}\def\Q{\mathbb Q} \def\texto{\text{o}} \def\LL{\Lambda}
\def\Rank{{\rm Rank}} \def\B{\scr B} \def\i{{\rm i}} \def\HR{\hat{\R}^d}
\def\to{\rightarrow}\def\l{\ell}\def\iint{\int}
\def\EE{\scr E}\def\Cut{{\rm Cut}}
\def\A{\scr A} \def\Lip{{\rm Lip}}
\def\BB{\scr B}\def\Ent{{\rm Ent}}

\maketitle

\begin{abstract} Let $P_t$ be  the (Neumann) diffusion semigroup $P_t$  generated by a weighted Laplacian on a complete connected Riemannian manifold $M$ without boundary or with a convex boundary. It is well known that the Bakry-Emery curvature is bounded below by a positive constant $\ll>0$ if and only if
$$W_p(\mu_1P_t, \mu_2P_t)\le \e^{-\ll t} W_p (\mu_1,\mu_2),\ \ t\ge 0, p\ge 1 $$ holds for all probability measures $\mu_1$ and $\mu_2$ on $M$,  where $W_p$ is the $L^p$ Wasserstein distance induced by the Riemannian distance. In this paper, we  prove the exponential contraction $$W_p(\mu_1P_t,  \mu_2P_t)\le c\e^{-\ll t} W_p (\mu_1,\mu_2),\ \   p\ge 1, t\ge 0$$ for some constants $c,\ll>0$ for a class of diffusion semigroups with negative curvature where the constant $c$ is essentially larger than $1$. Similar results are derived for SDEs with multiplicative noise by using explicit conditions on the coefficients, which are new even for SDEs with additive noise.
\end{abstract} \noindent
 \vskip 2cm

\section{Introduction}

Let $M$ be a $d$-dimensional connected complete Riemannian manifold possibly with a convex boundary $\pp M$.
Let   $\rr$ be the Riemannian distance. Consider    $L=\DD+Z$ for the Laplace-Beltrami operator $\DD$ and  some $C^1$-vector field $Z$ such that the (reflecting) diffusion process generated by $L$ is non-explosive.  Then the associated Markov semigroup $P_t$ is the (Neumann if $\pp M\ne\emptyset$) semigroup generated by $L$ on $M$. In particular, it is the case when  the curvature of $L$ is bounded below; that is,
\beq\label{C}\Ric_Z:= \Ric -\nn Z\ge K\end{equation} holds for some constant $K\in\R$. Here and throughout the paper,  we write  $\scr T\ge h$   for a (not necessarily symmetric)  $2$-tensor $\scr T$ and a function $h$ provided
$$\scr T  (X, X)  \ge h(x)|X|^2,\ \ X\in T_x M, x\in M.$$
There exist  many   inequalities on $P_t$  which are equivalent to the curvature condition \eqref{C},   see   \cite{BGL, O, RS, WBook} for details. In particular,  for any constant $K\in \R$, the
Wasserstein distance inequality
\beq\label{W} W_p ( \mu_1P_t, \mu_2P_t)\le \e^{-Kt} W_p (\mu_1,\mu_2),\ \ t\ge 0, p\ge 1, \mu_1,\mu_2\in \scr P(M)\end{equation} is equivalent to the curvature condition \eqref{C}.
  Here, $\scr P(M)$ is the class of all probability measures
on $M$;
 $W_p$ is the $L^p$-Warsserstein distance induced by $\rr$, i.e.,
$$W_p(\mu_1,\mu_2) :=\inf_{\pi\in \scr C(\mu_1,\mu_2)} \|\rr\|_{L^p(\pi)}, \ \ \mu_1,\mu_2\in \scr P(M), $$
where $\scr C (\mu_1,\mu_2)$  is the class of all couplings of $\mu_1$ and $\mu_2$;
and for   a Markov operator
 $P$ on $\B_b(M)$ (i.e. $P$ is a positivity-preserving linear operator with $P1=1$),
$$(\nu P)(A) := \nu(P 1_A),\ \ A\in \B(M), \nu\in \scr P(M),$$ where $\nu(f):=\int_M f\d\nu$ for $f\in L^1(\nu)$. In some references,
$\nu P$ is also denoted by $P^*\nu$. In the sequel we will use $P_t^*$ to stand for the adjoint operator of $P_t$ in $L^2(\mu)$ for the invariant probability measure $\mu$, hence adopt the notation $\nu P$ rather than $P^*\nu$ to avoid confusion.
When the curvature is positive (i.e. $K>0$), \eqref{W} implies the $W_p$-exponential contraction of $P_t$   for  $p\ge 1.$

\

In this paper, we aim to    consider the case when   \eqref{C} only holds for some negative constant $K,$ and to prove the exponential contraction
\beq\label{W'}W_p (\mu_1P_t,  \mu_2P_t)\le c\e^{-\ll t} W_p (\mu_1,\mu_2),\ \ t\ge 0, p\ge 1, \mu_1,\mu_2\in \scr P(M)\end{equation} for some constants $c,\ll>0$. It is crucial that the exponential rate $\ll$ is independent of $p$. Due to the equivalence of \eqref{C} and   \eqref{W}, in the negative curvature case it is essential that $c>1$.

According to \cite{WAnn},   even when $\Ric_Z$ is unbounded below,  i.e. $\Ric_Z$ goes to $-\infty$ when $\rr_o:=\rr(o,\cdot)\to\infty$ for a fixed $o\in M$, there may exist the log-Sobolev inequality which implies the exponentially convergence of $P_t$ in entropy. This suggests that \eqref{W'} may also hold for a class of diffusion semigroups with negative curvature.

Recently, some efforts have been made in this direction for $M=\R^d$, see   \cite{EB, EB2, LW}. More precisely, let $P_t$ be the diffusion semigroup for the solution to the following SDE on $\R^d$:
$$\d X_t = \ss 2\, \d B_t +b(X_t)\d t, $$where $B_t$ is the $d$-dimensional Brownian motion and   $b: \R^d\to \R^d$ is continuous.
If there exist constants $K_1,K_2, r_0>0$ such that
\beq\label{EB} \<b(x)-b(y), x-y\>\le  1_{|x-y|\le r_0} (K_1+K_2) |x-y|^2 - K_2 |x-y|^2,\ \ x,y\in \R^d,\end{equation}
then  due to \cite{EB,EB2} we have
\beq\label{Eb2} W_1(\dd_x P_t,  \dd_yP_t)\le c \e^{-\ll t} |x-y|,\ \ x,y\in \R^d, t\ge 0\end{equation}   for some constants $c,\ll>0$, where
$\dd_x$ is the Dirac measure at point $x$. Indeed, \cite{EB,EB2} proved the $W_1$-exponential contraction with respect to a modified distance
$f(|x-y|)$ in place of $|x-y|$ as constructed in \cite{CW, CW2} for estimates of the spectral gap using the coupling by reflection. Under condition \eqref{EB} the modified distance is comparable with the usual one so that \eqref{Eb2} follows. As mentioned in \cite{EB2} that  there is essential difficulty to prove \eqref{W'} for $p>1$  even for this flat case.

In Luo and Wang \cite{LW}  the estimate \eqref{Eb2} was extended as
\beq\label{Eb3} W_p( \dd_x P_t, \dd_yP_t)\le c \e^{-\ll t/p} (|x-y|+ |x-y|^{\ff 1 p}),\ \ x,y\in \R^d, t\ge 0, p\ge 1 \end{equation} for some constants $c,\ll>0$. Comparing with \eqref{W'} which is equivalent to
  $$W_p(\dd_x P_t, \dd_y P_t)\le c \e^{-\ll t} |x-y|,\ \ p\ge 1, x,y\in \R^d, t\ge 0$$   according to \cite{KK2} (see Proposition \ref{T3.1} below),
\eqref{Eb3} is less sharp for small $|x-y|$ and/or large $p$.  It is open whether \eqref{EB}, or
in the Riemannian setting that $\Ric_Z$ is uniformly positive outside a compact domain, implies \eqref{W'} for some constants $c,\ll>0$.

\

 As in \cite{KK, KK2}, we will consider the Warsserstein distances induced by Young functions in the class
\beg{align*} \scr N:= \Big\{\Phi\in C^1([0,\infty); [0,\infty)): &\ \Phi' \ \text{is\ nonnegative\ and \ increasing},\\
  \Phi(0)=0, &\Phi(r)>0\ \text{for}\ r>0, \lim_{r\to\infty} \ff{\Phi(r)}r=\infty\Big\}.\end{align*} For any $\Phi\in \scr N$ and a measure $\nu$ on $M$, consider the gauge norm  in $L^\Phi(\nu):$    $$\|f\|_{L^\Phi(\nu)}:= \inf\Big\{r>0: \nu\big(\Phi(r^{-1}|f|)\big)\le 1\Big\},\ \ \inf\emptyset:=\infty.$$
In particular,  we have $\|f\|_{L^{\Phi_p} (\nu)}=\|f\|_{L^p(\nu)}$ for  $\Phi_p(r):=  r^p,\ p\in (1,\infty)$.  This is the reason why we do not take $\Phi_p(r) =\ff 1 p r^p$ in the characterization of  Legendre conjugates.  We extend the notion $\Phi_p$ to $p= 1,\infty$ by letting $\Phi_1(r)=r, \Phi_\infty=\lim_{p\to\infty}\Phi_p$ and
$\|f\|_{L^{\Phi_p}(\nu)}=\|f\|_{L^p(\nu)}$ for all $p\in [1,\infty].$   Now, let $$W_\Phi(\mu_1,\mu_2)= \inf_{\pi\in \scr C(\mu_1,\mu_2)} \|\rr\|_{L^\Phi(\pi)},\ \
 \Phi\in \bar{\scr N}:= \scr N\cup\{\Phi_1,\Phi_\infty\}.$$ In particular, $W_{\Phi_p}=W_p$   for $p\in [1,\infty].$
We aim to prove the exponential decay
\beq\label{LL}W_\Phi( \dd_x P_t,  \dd_yP_t) \le c \Phi^{-1}(1)\e^{-\ll t} \rr(x,y),\ \ \ x,y\in M, t\ge 0, \Phi\in  \bar{\scr N}\end{equation} when   \eqref{C} only holds for a negative constant  $K$, where $\Phi^{-1}$ is the inverse of $\Phi (\ne \Phi_\infty)$ and  we set $\Phi^{-1}_\infty(1)=1$ by convention.

To extend condition \eqref{EB} to the Riemannian setting, consider the index
$$I(x,y)= \int_0^{\rr(x,y)} \sum_{i=1}^{d-1}\Big\{|\nn_{\dot\gg}J_i|^2 -\<\scr R(\dot\gg, J_i) \dot\gg,J_i\>\Big\}
 (\gg_s) \d s,\ \ x,y\in M,$$ where $\rr$ is the Riemannian distance, $\scr R$ is the curvature tensor;
$\gg: [0,\rr(x,y)]\to M$ is the minimal geodesic from $x$ to $y$ with unit speed;
    $\{J_i\}_{i=1}^{d-1}$ are   Jacobi fields along
$\gg$ such that $$J_i(y)=P_{x, y}\,J_i(x), \quad i=1,\ldots,d-1 $$ holds for the parallel transform
$P_{x,y}: T_x M\to T_y M$  along the geodesic $\gg$, and    $\{\dot \gg(s), J_i(s): 1\le i\le d-1\}$ ($s=0,\rr(x,y)$)  is an orthonormal basis of
the tangent space (at points $x$ and $y$, respectively).

Note that when $(x,y)\in \Cut(M)$, i.e. $x$ is in the cut-locus of $y$, the minimal geodesic may be not unique. As a convention  in the literature,   all conditions on the index $I$ are given outside $\Cut(M)$.
We now extend condition \eqref{EB} to the non-flat case as follows: for some constants $K_1,K_2>0$,
\beq\label{EB'} \beg{split}   I_Z(x,y)& := I(x,y)+ \<Z,\nn \rr(\cdot, y)\>(x)+ \<Z,\nn \rr(x,\cdot,)\>(y)\\
& \le \big\{(K_1+K_2)1_{\{\rr(x,y)\le r_0\}}- K_2\big\}\rr(x,y),\ \ x,y\in M.\end{split}\end{equation}   In the flat case we have $I(x,y)=0$ and $\rr(x,y)=|x-y|$,
  so that this condition reduces back to \eqref{EB}. Moreover, the curvature condition \eqref{C} is equivalent to
  $$I_Z(x,y)\le -K\rr(x,y),\ \ x,y\in M,$$ so that \eqref{EB'} implies $\Ric_Z\ge -(K_1+K_2).$

\

In the next section, we state our main results and present examples. With condition \eqref{EB'} we first   extend  the main results of \cite{EB, LW}  to  the present Riemannian setting and give  the exponential convergence of $P_t$ in $W_2$.   Under the ultracontractivity and  condition \eqref{C} for some $K<0$, our the second result ensures the desired inequality \eqref{LL}. Finally, we extend these results to SDEs with multiplicative noise by using explicit conditions on the coefficients. To prove these results, we make  some preparations   in Section 3. Complete proofs of the main results are addressed in Sections 4-6 respectively.

\section{Main Results and examples}

 We first consider the Riemannian setting, then extend to SDEs with multiplicative noise by using explicit conditions on the coefficients instead of  the less explicit curvature condition.

\subsection{The Riemannian setting}

We start with  condition \eqref{EB'}. Besides the extension of \eqref{Eb3},  this condition also implies  the hypercontractivity and the exponential convergence in $W_2$ for the semigroup $P_t$. For a measure $\mu$ and constants $p, q\ge 1$, let
 $\|\cdot\|_{L^p(\mu)\to L^q(\mu)} $ stand  for  the operator norm form $L^p(\mu)$ to $L^q(\mu)$. Recall that  $P_t$  is called hypercontractive
if it has a unique invariant probability measure $\mu$ and  $\|P_t\|_{L^2(\mu)\to L^4(\mu)}=1$ holds for large $t>0$. By interpolation theorem, $\|P_t\|_{L^2(\mu)\to L^4(\mu)}=1$ can be replaced by  $\|P_t\|_{L^p(\mu)\to L^q(\mu)}=1$ for some $\infty >q>p>1.$

\beg{thm}\label{T1.1}  Let  $\eqref{EB'}$ hold for some constants $K_1, K_2$ and $ r_0> 0$. Then:
\beg{enumerate} \item[$(1)$] There exist two constants $c,\ll>0$ such that for any $\Phi\in \bar{\scr N}$ and $x,y\in M$,
\beq\label{EB1'} W_\Phi(\dd_x P_t,  \dd_yP_t)\le  \inf\Big\{r>0:  \sup_{s\in (0, 1+ \rr(x,y)]} \ff{\Phi(r^{-1} s)}{s}\le \ff   {\e^{\ll t}}{c\rr(x,y)}\Big\},\ \ t\ge 0.\end{equation} In particular,
$$W_p(\dd_xP_t,  \dd_yP_t)\le \{c\e^{-\ll t}\}^{\ff 1 p} (\rr(x,y)+\rr(x,y)^{\ff 1 p}),\ \ p\ge 1, t\ge 0, x,y\in M.$$
\item[$(2)$] $P_t$ has a unique invariant probability measure $\mu$ and the log-Sobolev inequality
\beq\label{LS} \mu(f^2\log f^2) \le C \mu(|\nn f|^2) +\mu(f^2)\log \mu(f^2),\ \ f\in C_b^1(M)\end{equation} holds for some constant $C>0$. Consequently,  $P_t$  is hypercontractive.
\item[$(3)$] There exist  constants $c,\ll>0$ such that
\beq\label{W2} W_2( \nu P_t, \mu)\le c\e^{-\ll t} W_2(\nu, \mu),\ \ t\ge 0, \nu\in \scr P(M).\end{equation} \end{enumerate}
\end{thm}

To illustrate this result,   we present  below a consequence with   explicit curvature conditions in the spirit of \cite{WAnn}. These conditions allow    $\Ric_Z$ to be negative everywhere, for instance,
when $-C_1\le \Ric\le -C_2$  and  $C_2>-\nn Z\ge \dd$  for some constants $C_1>C_2>\dd>0$.
  As indicated in Introduction that  \eqref{EB'} implies  $\Ric_Z\ge -(K_1+K_2),$  so   in the following corollary we assume that $\Ric_Z$ is bounded below.

\beg{cor}\label{C1.2} Assume that   $\Ric_Z$ is bounded below.  Let   $\rr_o=\rr(o,\cdot)$ for a fixed point $o\in M$. If there exist constants $\si>0$ and  $\dd>\si (1+\ss 2)  \ss{d-1}$ such that
\beq\label{C1} -\nn Z \ge -\dd \ {\rm and} \  \Ric\ge  -\si^2\rr_o^2 \   {\rm   outside\ a\ compact\ set},\end{equation}  then
all assertions in Theorem $\ref{T1.1}$  hold.
\end{cor}

Next, we introduce sufficient conditions for \eqref{LL} which allow $\Ric_Z$ to be negative. Due to technical reason, we will need the ultracontractivity of $P_t$,  which is essentially stronger than the hypercontractivity.  We call $P_t$  ultracontractive if $\|P_t\|_{L^1(\mu)\to L^\infty(\mu)}<\infty$ for all $t>0.$ The ultracontractivity implies that $P_t$ has a density $p_t(x,y)$  with respect to $\mu$ (called heat kernel) and $$\|p_t\|_{L^\infty(\mu\times \mu)} = \|P_t\|_{L^1(\mu)\to L^\infty(\mu)}<\infty,\ \ t>0.$$ In references (see e.g. \cite{DS}), the ultracontractivity is also defined by $\|P_t\|_{L^2(\mu)\to L^\infty(\mu)}<\infty$ for $t>0$. When $P_t$ is symmetric in $L^2(\mu)$ we have
\beq\label{SY}\|P_t\|_{L^1(\mu)\to L^\infty(\mu)}\le\|P_{t/2}\|_{L^2(\mu)\to L^\infty(\mu)}^2,\ \ t>0,\end{equation} so that these two definitions  are equivalent. However, when $P_t$ is non-symmetric, the former might be stronger than the latter.  The appearance of the ultracontractivity  in our study is very nature: by Theorem \ref{T1.2}(1) we already have \eqref{LL} for $\Phi=\Phi_1$ (the weakest case), and by the ultracontractivity we are able to deduce the inequality from $\Phi_1$ to   $\Phi_\infty$ (the strongest case).    On the other hand, the result also indicates that \eqref{LL} implies the hypercontractivity of $P_t$.

\beg{thm}\label{T1.2}  Assume that   $\Ric_Z$ is bounded below. \beg{enumerate} \item[$(1)$] If    $P_t$ is ultracontractive, then
there exist constants $c,\ll>0$ such that for any $\Phi\in  \bar{\scr N} $,
\beq\label{LL0} W_{\Phi}( \dd_xP_t,  \dd_yP_t) \le \ff{c}{\Phi^{-1}(1)}  \e^{-\ll t} \min\Big\{ \rr(x,y),\  G_\Phi(t) \Big\},\ \ t>0, x,y\in M  \end{equation} holds for
$$G_\Phi(t):=\inf\Big\{r>0:\ (\mu\times\mu)\Big(\Phi\big(r^{-1}\rr\big)\Big)\le \|P_{t/2}\|^{-2}_{L^1(\mu)\to L^\infty(\mu)}\Big\}.$$ Consequently, for any $p\in [1, \infty], t\ge 0$ and $\mu_1,\mu_2\in \scr P(M)$,
\beq\label{LL1} W_p(\mu_1P_t,  \mu_2P_t)\le c \e^{-\ll t} \min\Big\{W_p(\mu_1,\mu_2),\  \|\rr\|_{L^p(\mu\times \mu)}\|P_{t/2}\|^{\ff 2p}_{L^1(\mu)\to L^\infty(\mu)}\Big\}.\end{equation}
\item[$(2)$] On the other hand, if there exist constants $c,\ll>0$ such that
\beq\label{WIF} W_\infty(\dd_xP_t, \dd_y P_t)\le c\e^{-\ll t} \rr(x,y),\ \ x,y\in M, t\ge 0,\end{equation} then the log-Sobolev inequality
$\eqref{LS}$ holds for $c=\ff{2c^2}\ll$, so that $P_t$ is hypercontractive.
\end{enumerate} \end{thm}

 We note that in  Theorem \ref{T1.2}(1) we have   $\|\rr\|_{L^p(\mu\times \mu)}<\infty$ for $p\in [1,\infty)$. Indeed,
 since $\Ric_Z$ is bounded below, by \cite[Theorem 2.1]{RW} the ultracontractivity implies the super log-Sobolev inequality \eqref{SLS} below, so that due to Herbst we have
$(\mu\times \mu)(\e^{r \rr^2})<\infty$ for all $r>0$ (see e.g. \cite{AMS}). Therefore, $G_\Phi(t)<\infty$ for $t>0$ and $\Phi\in \scr N$ satisfying
$$\limsup_{r\to\infty}\ff{\log \Phi(r)}{r^2}<\infty.$$

In the symmetric case (i.e. $Z=\nn V$ for some $V\in C^2(M)$), explicit sufficient conditions  for the ultracontractivity have been introduced
in   \cite{WAnn} by using  the dimension-free Harnack inequality in the sense of \cite{W97}.  Together with a suitable  exponential  estimate on the diffusion process, this inequality implies $\|P_t\|_{L^2(\mu)\to L^\infty(\mu)}<\infty$ for $t>0$  and thus,   $P_t$ is ultracontractive   due to \eqref{SY}.
The conditions can be formulated as \beq\label{C2} -\nn Z\ge \Psi_1\circ\rr_o\ {\rm and}\ \Ric \ge
-\Psi_2\circ\rr_o\ \text{hold \ outside\ a\
compact\ subset\ of\ }M,\end{equation}
where $\Psi_1, \Psi_2: (0,\infty)\to (0,\infty)$ are increasing functions such that
\beq\label{4.3}\int_1^\infty \ff{\d s}{\ss s\int_0^{\ss
s}\Psi_1(u)\d u}<\infty,\ \  \lim_{r\to\infty}\min\Big\{\Psi_1(r), \ff{(\int_0^r \Psi_1(s)\d s)^2}{\Psi_1(r)}\Big\} =\infty,\end{equation}
and for some constants $\theta\in (0, 1/(1+\ss 2))$ and $C>0,$
\beq\label{4.4} \ss{\Psi_2(r+t)(d-1)}\le \theta \int_0^r \Psi_1(s)\d s
+\ff 1 2 \int_0^{t/2}\Psi_1(s)\d s+C,\ \ r,t\ge 0.\end{equation}  When $\Ric$ is bounded below,  \eqref{4.4} as well as  the second inequality in \eqref{C2} hold for  $\Psi_2$ being a large enough constant. In general, since $\int_0^r \Psi_1(s)\d s \ge 2 \int_0^{r/2} \Psi_1(s)\d s$, \eqref{4.4} with $\theta=\ff 1 4<\ff 1 {1+\ss 2}$ follows from
 \beq\label{4.4'}\beg{split}  \ss{\Psi_2(r)(d-1)}&\le \ff 1 2\inf_{t\in [0,r]} \bigg\{ \int_0^{t/2} \Psi_1(s)\d s
+  \int_0^{(r-t)/2}\Psi_1(s)\d s\bigg\}+C\\
&= \int_0^{r/4} \Psi_1(s)\d s+C,\ \ r\ge 0.\end{split}\end{equation}

Since   \eqref{SY} fails for   non-symmetric semigroups, we apply the inequality
$$\|P_t\|_{L^1(\mu)\to L^\infty(\mu)}\le  \|P_{t/2}\|_{L^1(\mu)\to L^2(\mu)}\|P_{t/2}\|_{L^2(\mu)\to L^\infty(\mu)}$$ due to the semigroup property. So,  to ensure  the ultracontractivity,  we need an additional condition  implying    $\|P_t\|_{L^1(\mu)\to L^2(\mu)}<\infty$  (see Corollary \ref{C1.3}(2) below).

 To estimate $G_\Phi(t)$ in \eqref{LL0} using $\Psi_1$, we introduce
$$\LL_1(r):= \ff 1 {\ss r} \int_0^{\ss r} \Psi_1(s)\d s,\ \ \
 \LL_2(r):= \int_r^\infty \ff{\d\ s}{\ss s\int_0^{\ss s}\Psi_1(u)\d
u},\ \ \ r>0.$$ Obviously,  the inverse function $\LL_2^{-1}$ exists on $(0,\infty)$, and since $\LL_1$ is increasing with $\LL_1(\infty)=\infty$, we have
 $$\LL_1^{-1}(r):=\inf\{s\ge 0: \LL_1(s)\ge r\}<\infty,\ \ r\ge 0.$$

\beg{cor}\label{C1.3} Assume that     $\eqref{4.3}$ and $\eqref{4.4}$ hold for some constants $\theta\in (0, 1/(1+\ss 2))$ and $C>0.$
\beg{enumerate} \item[$(1)$] If $P_t$ is symmetric, i.e.   $Z=\nn V$ for some $V\in C^2(M)$, then
  there exist constants $c,\ll>0$ such that $\eqref{LL0}$ and $\eqref{LL1}$    hold  for
$$G_\Phi(t):=\inf\Big\{\ll>0:\ (\mu\times\mu)\big(\Phi (\ll^{-1}\rr  )\big)\le \e^{-  ct^{-1}\{1+\LL_1^{-1}(c t^{-1})-\LL_2^{-1} ( c^{-1}t  )\}}\Big\},\ \ t>0.$$
 \item[$(2)$] If  $P_t$ is non-symmetric but there exists continuous $h\in C([0,1]; [0,\infty)) $ with $h(r)>0$ for $r>0$ such that
$\int_0^1 \ff{ h(r)} r\d r<\infty$ and
$$H(\theta):=\int_0^1 \ff \theta {h(r)}\Big\{1+\LL_1^{-1} \big(  \theta/ h(r)\big) +  \LL_2^{-1} \big( h(r)/\theta\big)\Big\}\d r<\infty,\ \ \theta>0,$$ then there exist constants $c,\ll>0$ such that $\eqref{LL0}$ holds for
$$G_\Phi(t):=\inf\Big\{\ll>0:\ (\mu\times\mu)\big(\Phi (\ll^{-1}\rr  )\big)\le \e^{-  ct^{-1}\{1+\LL_1^{-1}(c t^{-1})-\LL_2^{-1} ( c^{-1}t  )\}-c H(ct^{-1})}\Big\}.$$
 \end{enumerate}
\end{cor}

To conclude this part, we
  present a simple example to illustrate Corollary \ref{C1.3}.

\paragraph{Example 2.1.} Let $M$ have   non-positive sectional curvatures and a pole $o\in M$. Let $Z= Z_0-\dd \nn \rr_o^{2+\vv}$ outside a compact domain, where   $\dd,\vv>0$ are constants and   $Z_0$ is a $C^1$ vector field with
\beq\label{Z1}\limsup_{\rr_o\to\infty}\ff{|\nn Z_0|}{\rr_o^\vv}<\dd(1+\vv)(2+\vv).\end{equation}  Let $\Psi_2: (0,\infty)\to (0,\infty)$ be increasing such that
\beq\label{Z2} \Ric \ge - \Psi_2(\rr_o), \ \lim_{r\to\infty} \ff{\Psi_2(r)}{r^{2(1+\vv)}}=0. \end{equation}
By  \eqref{Z1}, \eqref{Z2} and the Hessian comparison theorem, we see that   \eqref{C2}, \eqref{4.3} and \eqref{4.4'} hold with
$\Psi_1(r)= c_1 r^\vv$ for some constant $c_1>0$.
According to Corollary \ref{C1.3},  there exist constants $c,\ll>0$ such that for any $p\ge 1$,
$$W_p(\mu_1P_t,  \mu_2P_t) \le c \e^{-\ll t} \min\Big\{W_p(\mu_1,\mu_2),\  \|\rr\|_{L_p(\mu\times\mu)}\exp\Big[\ff {c  } {pt^{1+\ff 2\vv}}\Big] \Big\},\ \ t>0, \mu_1,\mu_2\in \scr P(M).$$

\subsection{SDEs with multiplicative noise}

Consider the following SDE on $\R^d$:
\beq\label{SDE} \d X_t= b(X_t)\d t +\ss 2 \si(X_t)\d B_t,\end{equation} where $B_t$ is the $m$-dimensional Brownian motion, $b: \R^d\to \R^d$ and $\si: \R^d \to  \R^d\otimes\R^m $ (the space of $d\times m$-matrices) are locally Lipshitz such that
$$\|\si\|_{HS}^2(x) +  \<b(x),x\>\le C(1+|x|^2),\ \ x\in \R^d$$ holds for some constant $C>0$, where and in the following, $\|\cdot\|_{HS}$ and $\|\cdot\|$ denote  the Hilbert-Schmidt and the operator norms respectively.  Then the SDE has a unique solution $\{X_t(x)\}_{t\ge 0}$ for every initial point $x\in \R^d$. Let $P_t$ be the associated Markov semigroup:
$$P_t f(x):= \E[f(X_t(x))],\ \ t\ge 0, x\in \R^d, f\in \B_b(\R^d).$$
We intend to investigate the $W_p$-exponential contraction for $p\in [1,\infty)$. As mentioned in Introduction that existing results only apply to $p=1$ and $\si=I$, and as mentioned in \cite{EB2,LW} that there is essential difficulty to prove \eqref{W'} for $p>1$ even for $\si=I$.   So, the
present study is non-trivial.

Corresponding to that \eqref{C} implies \eqref{W} in the Riemannian setting, we have the following assertion.

\beg{thm}\label{ST1} Let $p\in [1,\infty)$. If
\beq\label{DSS} \beg{split} & \ff{(p-2) |(\si(x)-\si(y))^*(x-y)|^2}{|x-y|^2}+ \|\si(x)-\si(y)\|_{HS}^2 +\<b(x)-b(y), x-y\>\\
&\le -K_p|x-y|^2,\ \ x\ne y\in \R^d\end{split}\end{equation} holds for some constant $K_p\in \R$, then
$$ W_{p}(\mu_1P_t,  \mu_2P_t)\le \e^{-K_p t} W_{p}(\mu_1,\mu_2),\ \ \ t\ge 0, \mu_1,\mu_2\in \scr P(\R^d).$$
 \end{thm}
Note that this result does apply to $p=\infty$ when $\si$ is non-constant.
Next, as in the Riemannian case, we intend to prove the exponential contraction in $W_p$ when \eqref{DSS} only holds  for some negative constant $K_p$.  To this end, we need the SDE to be non-degenerate.   The following result contains   analogous assertions in  Theorems \ref{T1.1} and \ref{T1.2},  where the first assertion extends \eqref{Eb2} to the multiplicative noise setting.

\beg{thm}\label{ST2}  Assume that
$ \si\si^* \ge \ll_0^2 I$  for some constant $\ll_0>0$.     \beg{enumerate} \item[$(1)$] If there exist constants $K_1,K_2,r_0>0$ such that  $Z$ and
$\si_0:= \ss{\si\si^* -\ll_0^2 I}$ satisfy
\beq\label{DSS2} \beg{split} & \|\si_0(x)-\si_0(y)\|_{HS}^2-\ff{|(\si(x)-\si(y))^*(x-y)|^2}{|x-y|^2} +\<b(x)-b(y), x-y\>\\
&\le \big\{(K_1+K_2)1_{\{|x-y|\le r_0\}} - K_2\big\} |x-y|^2,\ \ x,y\in \R^d,\end{split}\end{equation}  then there exist constants $c,\ll>0$ such that
$$ W_1(\mu_1 P_t, \mu_2P_t)\le c \e^{-\ll t} W_1(\mu_1,\mu_2),\ \ \ t\ge 0, \mu_1,\mu_2\in \scr P(\R^d).$$
\item[$(2)$] Let $P_t$ have a unique invariant probability measure $\mu$ such that the log-Sobolev inequality
\beq\label{LSA} \mu(f^2\log f^2)\le C\mu(|\si^*\nn f|^2),\ \ f\in C_b^1(\R^d), \mu(f^2)=1\end{equation} holds for some constant $C>0$.
 If there exists a constant $K>0$ such that
 \beq\label{DSS2'}   \|\si(x)-\si(y)\|_{HS}^2 +\<b(x)-b(y), x-y\>
 \le  K |x-y|^2,\ \ x,y\in \R^d, \end{equation}    then $\eqref{W2}$ holds for some   constants $c,\ll>0$ and $M=\R^d$.
 \item[$(3)$] Let  $P_t$ be ultracontractive and let $\eqref{DSS2'}$ hold for some constant $K>0$. Then there exist a constant $\ll>0$ such that for any $p\in [1,\infty)$, condition $\eqref{DSS}$  implies   $\eqref{LL1}$ for   some constant  $c=c(p)>0$, and all $t\ge 0, \mu_1,\mu_2\in \scr P(\R^d)$.
\end{enumerate}
\end{thm}
According to \cite[Lemma 3.3]{PW}, we have
\beq\label{EP} \|\si_0(x)-\si_0(y)\|^2\le \ff 1 {4\ll_0} \|(\si\si^*)(x)-(\si\si^*)(y)\|_{HS}^2,\ \ x,y\in \R^d.\end{equation}
Combining this with $\|\cdot\|_{HS}^2 \le d \|\cdot\|^2$, we see that  \eqref{DSS2} follows from the following more explicit condition:
\beq\label{DSS3} \beg{split} &\ff {d-1} {4\ll_0} \|(\si\si^*)(x)-(\si\si^*)(y)\|_{HS}^2  +\<b(x)-b(y), x-y\> \\
& \le \big\{(K_1+K_2)1_{\{|x-y|\le r_0\}} - K_2\big\} |x-y|^2,\ \ x,y\in \R^d. \end{split}\end{equation}

Note that conditions in Theorem \ref{ST1} and Theorem \ref{ST2}(1) are explicit. To illustrate Theorem \ref{ST2}(2)-(3),   we present below sufficient conditions for the log-Sobolev inequality \eqref{LSA} and the ultracontractivity of $P_t$. For $a:=\si\si^*$ and $(g_{ij})_{1\le i,j\le d}:= a^{-1}$, we introduce the   Christoffel symbols
$$\GG_{ij}^k:= \ff 1 2 \sum_{m=1}^d \big(\pp_i g_{mj}+\pp_j g_{im} - \pp_m g_{ij}\big)a_{km},\ \ 1\le i, j, k\le d,$$ and the matrix $\GG ab$:
$$(\GG a b)_{ij}:= \sum_{k,l=1}^d \GG_{kl}^i a_{kj} b_k,\ \ 1\le i,j\le d.$$

\beg{prp}\label{PPN} Let $\si\in C_b^2(\R^d\to \R^d\otimes \R^d)$ such that $a:=  \si\si^*\ge \aa I$ for some constant $\aa>0$, and let $b\in C^1(\R^d\to \R^d)$ such that
\beq\label{CC0}  \ff 1 2  (\GG a b +\nn_b a) -   (\nn b)a\ge -K_0 I \end{equation} for some constant $K_0$.   If there exist  constants $c_1,c_2>0$ and $\dd>1$   such that
\beq\label{UL} L |\cdot|^2 \le c_1 -c_2|\cdot|^{2\dd}, \end{equation}
then $P_t$ has a unique invariant probability measure $\mu$ and there exists a constant $c>0$ such that $$\|P_t\|_{L^1(\mu)\to L^\infty(\mu)}\le \exp\Big[c+ct^{-\ff{\dd}{\dd-1}}\Big],\ \ t>0.$$
\end{prp}

  We now introduce a simple example to illustrate Theorem \ref{ST2}.

\paragraph{Example 2.2.}  Let $\si\in C_b^2(\R^d\to \R^d\otimes \R^d)$ such that $a:=  \si\si^*\ge \aa I$ for some constant $\aa>0$.
Let $b(x)= -c_0 |x|^\theta x $  for   large $|x|,$  where $c_0>0$ and $\theta>0$ are constants.  Obviously, condition \eqref{DSS2'} holds. If
\beq\label{*D} \lim_{|x|\to\infty} |x|\cdot\|\nn \si(x)\| =0,\end{equation}  then \eqref{CC0} holds for some constant $K_0$.  Moreover, it is easy to see that
$$L  |\cdot|^2  \le  c_1-c_2 |x|^{\theta+2},\ \ \ll>0, x\in\R^d$$ holds for some constants $c_1,c_2>0$.
By Proposition \ref{PPN}  and Theorem \ref{ST2}(3), for any $p\in [1,\infty)$, there exist constants $\ll,c>0$ such that
$$W_p(\mu_1P_t,  \mu_2P_t)\le c  \e^{-\ll t} \min\Big\{W_p(\mu_1,\mu_2),\ \exp\Big[ c  t^{-\ff{\theta+2}\theta}\Big] \Big\},\ \ t> 0, \mu_1,\mu_2\in \scr P(\R^d).$$

\section{Preparations}

This section includes some propositions which will be used to prove the   results introduced in Section 2. We first  recall a link between  the Wasserstein distance and gradient estimates due to \cite{KK2},  then  deduce the hyperboundedness and the exponential convergence in entropy from the log-Sobolev inequality for non-symmetric diffusion semigroups, and finally prove the exponential contraction in gradient  for ultracontractive semigroups  in a general framework including both diffusion and jump Markov semigroups.

\subsection{Wasserstein distance and gradient inequalities}

Let $(E,\rr)$ be  a geodesic Polish space, i.e. it is a Polish space  and for any two different points $x,y\in E$, there exists a continuous curve $\gg: [0,1]\to E$ such that $\gg_0=x, \gg_1=y$ and $\rr(\gg_s,\gg_t)=|s-t|\rr(x,y)$ for $s, t\in [0,1].$ Then for any $f\in \Lip_b(E)$, the class of bounded Lipschitz functions on $E$, the length of gradient
$$|\nn f|(x):= \limsup_{\rr(x,y)\downarrow 0} \ff{|f(x)-f(y)|}{\rr(x,y)},\ \ x\in E$$ is measurable.
Moreover, let $P(x,\d y)$ be a Markov transition kernel and define the Markov operator
$$Pf(x):= \int_E f(y)P(x,\d y),\ \ x\in E, f\in \B_b(E).$$
For any $\Phi\in\bar{\scr N}\setminus\{\Phi_\infty\}$, consider the Young norm induced by $\Phi$ with respect to $P$
\beq\label{*NP}\|f\|_{L_*^\Phi(P)}(x) :=\sup\Big\{P(fg)(x): g\in \scr B_b(E), P\Phi(|g|)(x)\le 1\Big\},\ \ x\in E, f\in \B_b(E), \end{equation}
and  set $\|f\|_{L_*^{\Phi_\infty}(P)}(x)= P|f|(x).$ Then $\|\cdot\|_{L_*^{\Phi_p}}=\|\cdot\|_{L^{\Phi_q}}$ for $p\in [1,\infty], q=\ff{p}{p-1}.$
The following result follows from \cite[Theorem 2.2,  Remark 2 and Remark 3]{KK2}.

\beg{prp}[\cite{KK2}] \label{T3.1} For any constant $C>0$ and $\Phi\in \bar{\scr N}$,  the following statements are equivalent to each other:
\beg{enumerate} \item[$(1)$]
$ |\nn P f| \le C \|\nn f\|_{L_*^\Phi(P)} $ for $f\in \Lip_b(E).$
\item[$(2)$] $W_{\Phi}(\dd_xP,  \dd_yP) \le C \rr(x,y),\ \  x,y\in  E.$
\end{enumerate} When $\Phi=\Phi_p$ for $p\in [1,\infty]$, they are also equivalent to
\beg{enumerate} \item[$(3)$] $W_{p}(\mu_1P,  \mu_2P) \le C W_{p} (\mu_1,\mu_2),\ \ \mu_1,\mu_2\in \scr P(E).$\end{enumerate}
\end{prp}

\subsection{Hyperboundedness and exponential convergence in entropy}

 When $P_t$ is symmetric, it is well known that the hyperbounddeness, exponential convergence in entropy and   the log-Sobolev inequality are equivalent each other, see    \cite{BGL, Wbook} and references within. In   the non-symmetric case, the log-Sobolev inequality implies the former two properties  if the generator $L$ and the symmetric part of the Dirichlet form $\EE$ satisfy
\beq\label{NC} \beg{split} &-\mu((1+\log f)Lf)\ge c_0\EE\big(\ss f,\ss f\big) \ {\rm and}\\
& -\mu(f^{p-1} Lf)= \ff{c_0(p-1)}{p^2} \EE( f^{\ff p 2},  f^{\ff p 2}),\ \ p>1, f\in \D\end{split}\end{equation}  for some constant $c_0>0$ and a reasonable class  $\D$ of non-negative bounded functions, which is stable under $P_t$ and dense in $L^p_+(\mu)
:=\{f\in L^p(\mu): f\ge 0\}$ for any $p\ge 1$, see e.g. \cite{Gross}. In applications, it may be not easy to   figure out the class $\D$ such that \eqref{NC} holds. But in general this condition can be replaced by  the following approximation formula Lemma \ref{L0} in the spirit of \cite{RW04}.

Now, consider the (Neumann) semigroup $P_t$ generated by $L:= \DD+Z$ for a local bounded vector field $Z$ such that $P_t$ has a unique invariant probability measure $\mu$.
Let $$\D_0=\big\{f\in C_0^\infty(M):\ f\ \text{satisfies\ the\ Neumann\ condition\ if\ } \pp M\ne\emptyset\big\}.$$ Then $(L,\D_0)$ is dissipative
(thus, closable)  in $L^1(\mu)$ with closure $(L,\D_1(L))$ generating $P_t$ in $L^1(\mu)$, see e.g. \cite{ST} and references within. Let
$$\D= \{f\in \D_1(L)\cap L^\infty(\mu):\ f\ge 0\}. $$

\beg{lem}\label{L0} Let $f\in \D$ and $\psi \in C_b^\infty([{\rm ess}_\mu\inf f, \infty))$. There exists a sequence $\{f_n\}_{n\ge 1}\subset \D_0$ with $\inf f_n=\inf f$ such that $f_n\to f$ in $L^m(\mu)$ for any $m\ge 1$, $L f_n\to Lf$ in $L^1(\mu)$, and
$$\mu(\psi(f)Lf) = -\lim_{n\to\infty} \mu(\psi'(f_n)|\nn f_n|^2).$$\end{lem}
\beg{proof}  Since $f\in \D\subset \D_1(L)\cap L^\infty(\mu)$, there exists a uniformly bounded sequence $\{f_n\}_{n\ge 1}\subset \D_0$ such that $\inf f_n= {\rm ess}_\mu\inf f$ and $f_n\to f, L f_n\to Lf$ in $L^1(\mu)$.   By the uniform boundedness,  $f_n\to f$ in $L^m(\mu)$ for any $m\ge 1$. Since $\psi\in C_b^\infty([\inf f_n, \infty))$,
$$g_n:= \int_{\inf f_n}^{f_n}\psi(s)\d s\in \D_c:=\{g+c:\ c\in \R, g\in \D_0\}\subset \D_1(L).$$ This implies
$\mu(L g_n)=0$ since $\mu$ is $P_t$-invariant. So, by the dominated convergence theorem,
$$\mu(\psi(f)Lf)= \lim_{n\to\infty} \mu(\psi(f_n)L f_n) = \lim_{n\to\infty} \mu(L g_n - \psi'(f_n)|\nn f_n|^2) = - \lim_{n\to\infty} \mu(\psi'(f_n)|\nn f_n|^2).$$
\end{proof}

\beg{prp}\label{PN} Let $Z$ be a locally bounded vector field  such that the (Neumann) semigroup $P_t$ generated by  $L:=\DD+Z$  has  a unique invariant probability measure $\mu$.  \beg{enumerate} \item[$(1)$] If the super log-Sobolev inequality
\beq\label{SLS} \mu(f^2\log f^2)\le r \mu(|\nn f|^2) +\bb(r),\ \ r>0,\ \ f\in C_b^1(M), \mu(f^2)=1. \end{equation}
holds for some $\bb\in C((0,\infty);(0,\infty))$, then for any constants    $q>p\ge 1$ and   $\gg\in C((p,q); (0,\infty))$ such that
$t:=\int_p^q \ff{\gg(r)}r\d r<\infty,$  there holds
$$\|P_t\|_{L^p(\mu)\to L^q(\mu)}\le \exp\bigg[\int_p^q\ff{\bb(4\gg(r)(1-r^{-1}))}{r^2}\,\d r\bigg].$$
\item[$(2)$] If the log-Sobolev inequality
\beq\label{LS} \mu(f^2\log f^2)\le C \mu(|\nn f|^2) + \mu(f^2)\log \mu(f^2),\ \ f\in C_b^1(M) \end{equation} holds for some constant $C>0$,
then
$$\mu((P_t g)\log P_t g)\le \e^{-4t/C} \mu(g\log g),\ \ g\in \B_b(M), g\ge 0, \mu(g)=1.$$ \end{enumerate} \end{prp}

\beg{proof}  (1) According to Lemma \ref{L0}, for any $f\in \D$ and $p>1$, there exists $\{f_n\}_{n\ge 1}\subset \D_0$ such that
$f_n\to f^{\ff p 2}$ in $L^m(\mu)$ for all $m\ge 1$, and
\beq\label{IT} -\mu(f^{p-1} Lf)= \ff{4(p-1)}{p^2} \limsup_{n\to\infty} \mu(|\nn f_n|^2).\end{equation}
Applying \eqref{SLS} to $f_n$ and using \eqref{IT}, we obtain
\beg{align*} & p\mu(f^p\log f)=\lim_{n\to\infty} \mu(f_n^2\log f_n^2)
\le    r \liminf_{n\to\infty} \mu(|\nn f_n|^2)+\bb(r) \\
&\le  \ff{rp^2}{4(p-1)}\Big(-\mu(f^{p-1} Lf)+\ff{4\bb(r)(p-1)}{rp^2}\Big),\ \ r>0.\end{align*}
Set $c(p) =\ff{rp}{4(p-1)},$  we have
$$\ff{4\bb(r)(p-1)}{r p^2}= \ff{\bb(4c(p)(1-p^{-1}))}{pc(p)},\ \ p>1,$$ so that the above inequality becomes
$$\mu(f^p\log f) \le c(p) \Big(-\mu(f^{p-1} Lf) +\gg(p)\Big),\ \ p>1, f\in \D$$ for
$\gg(p):= \ff{\bb(4c(p)(1-p^{-1}))}{pc(p)}.$ Noting that $\D$ is $P_t$-invariant (i.e. $P_t\D\subset \D$)  and   dense in
$L_+^p(\mu)$ for any $p\ge 1$,  the desired  assertion follows from the proof of \cite[Corollary 3.13]{Gross}.

(2)  It suffices to prove for $g\in \D$ with $\inf g>0.$ Applying Lemma \ref{L0} to $f= P_tg$ and $\psi(s)= 1+ \log s$, and using  \eqref{LS}, we obtain
\beg{align*} &\ff{\d}{\d t} \mu( (P_t g)\log P_t g) = \mu((1+\log P_t g)  L P_t g) = -4\lim_{n\to\infty} \mu\big(\big|\nn\ss f_n\big|^2\big)\\
& \le   -\ff 4 C\liminf_{n\to\infty}  \mu(f_n\log f_n) = -\ff 4 C \mu((P_t g)\log P_t g),\ \ t\ge 0.\end{align*}
This implies the desired exponential estimate. \end{proof}

\subsection{Exponential contraction in gradient}

In this part, we consider a general framework including both diffusion and   jump processes.
Let $(E,\F,\mu)$ be a separable complete probability space, and let $P_t$ be a Markov semigroup on $L^2(\mu)$ with $\mu$ as invariant probability measure. Let $(L,\D(L))$ be the generator of $P_t$ in $L^2(\mu)$. We assume that there exists an algebra $\scr A\subset \D(L)$ such that
\beg{enumerate} \item[{\rm (i)}] $1\in \scr A$, $\scr A$   is dense in $L^2(\mu)$ and  the algebra induced by
$$\D:=\{P_s f: s\ge 0, f\in \scr A\}$$ is contained in $\D(L)$.
\item[{\rm (ii)}] $\GG(f,g):= \ff 1 2 (L(fg) - fLg -gLf)$
gives rise to a non-degenerate positive definite bilinear form on $\scr D\times\scr D$; i.e., for any $f\in \scr D$, $\GG(f,f)\ge 0$ and it equals to $0$ if and only if $f$ is constant.  \end{enumerate}
In particular, when $P_t$ is the (Neumann) semigroup generated by $L:=\DD+Z$ on $M$ with $\Ric_Z$ bounded below, the assumption holds for
$$\scr A:= \{f+c: f\in C_0^\infty(M) \ {\rm satisfying \ the\ Neumann\ condition\ if} \ \pp M\ne \emptyset, c\in\R\}.$$
Under the above conditions,
$$\EE(f,g):=\mu(\GG(f,g)),\ \ f,g\in \scr A$$ is closable and the closure $(\EE, \D(\EE))$ is a conservative symmetric Dirichlet form.
Although $P_t$ is not associated to $(\EE, \D(\EE))$ when it is non-symmetric,  we   have
\beq\label{REL} \ff{\d }{\d t} \mu((P_t f)^2) = - 2 \EE(P_tf, P_t f),\ \ t\ge 0, f\in \scr D.\end{equation}

 If $\|P_t\|_{L^1(\mu)\to L^\infty(\mu)}<\infty,$ then $P_t$ has a heat kernel
$p_t(x,y)$ with respect to $\mu$, i.e.
$$P_t f=\int_E p_t(\cdot, y)f(y)\mu(\d y),\ \ f\in L^2(\mu),$$ and
$${\rm ess}_{\mu\times \mu}\sup p_t=\|P_t\|_{L^1(\mu)\to L^\infty(\mu)}<\infty.$$
We consider the $``$gradient" length $|\nn_\GG f|= \ss{\GG(f,f)}$ induced by $\GG$. Note that for jump processes the length is non-local  and thus essentially different from the usual gradient length. As shown below that   estimates of $|\nn_\GG P_t|$ have a close  link to functional inequalities of the associated Dirichlet form.

\beg{prp}\label{P2.1} Assume that there exist $t_1>0$ and $\eta \in C([0,\infty); (0,\infty))$ such that
\beq\label{AS} \|P_{t_1}\|_{L^1(\mu)\to L^\infty(\mu)} <\infty,\ \ |\nn_\GG P_t f|^2\le \eta(t) P_t |\nn_\GG f|^2,\ t\ge 0, f\in \scr D.\end{equation}
Then there exist constants $c,\ll, t_2>0$ such that for any $q\ge 1$ and $\eta_q\in C([0,\infty); (0,\infty))$, the gradient estimate    \beq\label{GP} |\nn_\GG P_t f|^2\le \eta_q(t)  (P_t |\nn_\GG f|^{q})^{\ff 2 q},\ t\ge 0, f\in\D\end{equation} implies
\beq\label{GE} \|\nn_\GG P_tf \|_{L^\infty(\mu)}^2 \le \Big(c\sup_{[0,t_2]} \eta_q\Big) \e^{-\ll t} {\rm ess}_\mu\inf(P_t |\nn_\GG f|^{q})^{\ff 2 q},
\ \ t\ge t_2, f\in \D.\end{equation}
\end{prp}

\beg{proof} (a) We first prove
\beq\label{EP} \EE(P_tf, P_tf) \le C \e^{-\ll t} \EE(f,f),\ \ f\in \D, t\ge 0\end{equation} for some constants $C,\ll>0$.
By the second inequality in \eqref{AS},  for any $t>0$ and $f\in \D$  we have
 $$\ff{\d }{\d s} P_s (P_{t-s}f)^2 = 2 P_s |\nn_\GG P_{t-s} f|^2 \le 2 \eta(t-s)  P_t |\nn_\GG f|^2,\ \ s\in [0,t].$$ Integrating
both sides over $[0,t]$  leads to
$$ P_t f^2\le (P_t f)^2 +  C(t)P_t |\nn_\GG f|^2,\ \ C(t):=2\int_0^t \eta(s) \d s,\  t>0. $$   Taking   $t=t_1$ and noting that $\mu$ is the invariant probability measure of $P_t$, we obtain
\beq\label{DP} \mu(f^2)\le C(t_1) \EE(f,f)+ \|P_{t_1}\|_{1\to\infty}^2 \mu(|f|)^2,\ \ f\in \D.\end{equation}
Since $\D(\EE)$ is the closure of $\D$ under the $\EE_1$-norm, this inequality also holds for   $f\in \D(\EE).$  By condition (ii), the symmetric    Dirichlet form is irreducible. So, according to  \cite[Corollary 1.2]{W14} the defective Poincar\'e inequality \eqref{DP} implies
the Poincar\'e inequality
\beq\label{P} \mu(f^2)\le \ff 1 {\ll} \EE(f,f) +\mu(f)^2,\ \ f\in \D(\EE) \end{equation}  for some constant $\ll>0$.
By \eqref{REL} and that $\D$ is dense in $L^2(\mu)$, the Poincar\'e  inequality is equivalent to
\beq\label{EXP} \|P_t f-\mu(f)\|_2\le \e^{-\ll t} \|f-\mu(f)\|_2,\ \ t\ge 0, f\in L^2(\mu).\end{equation}
On the other hand, by the second inequality in \eqref{AS}, for any $t>0$ and $f\in \D$  we have
 $$\ff{\d }{\d s} P_s (P_{t-s}f)^2 = 2 P_s |\nn_\GG P_{t-s} f|^2 \ge \ff 2 {\eta(s)}  |\nn_\GG P_t f|^2,\ \ s\in [0,t].$$  So,
$$|\nn_\GG P_tf|^2 \le \ff {P_t f^2-(P_tf)^2} {2\int_0^t \eta(s)^{-1} \d s},\ \ t>0, f\in \D.$$
Using $P_tf-\mu(f)$ to replace $f$ and integrating  with respect to $\mu$, we obtain
$$\EE(P_{2t}f, P_{2t}f) \le \ff {\|P_tf-\mu(f)\|_2^2} {2\int_0^t \eta(s)^{-1} \d s},\ \ t>0, f\in \D.$$
Combining this with \eqref{EXP} and \eqref{P} we arrive at
$$ \EE(P_{t}f, P_{t}f) \le c_1 \e^{-\ll t} \EE(f,f),\ \ t\ge 1, f\in \D$$  for some constant $c_1>0$; that is,      \eqref{EP} holds for $t>1.$ Finally, \eqref{AS} implies \eqref{EP} for $t\in [0,1].$

(b) Next, we intend to find out a constant $t_0\ge t_1$ such that
\beq\label{H} \ff 1 2 \le p_t \le 2,\ \ (\mu\times\mu){\text -a.e.}, t\ge t_0.\end{equation}
Indeed, by \eqref{EXP} and the first inequality in \eqref{AS}, we obtain
\beg{align*} &\bigg|\int_E (p_{t+2 t_1}(\cdot, y)-1)f(y)\mu(\d y)\bigg|= |P_{t_1}(P_{t+t_1}f-\mu(f))| \\
&\le c_0 \mu(|P_{t+t_1}f-\mu(f)|)
 \le c_0 \e^{-\ll t} \|P_{t_1}f-\mu(f)\|_2\le c_0^2 \e^{-\ll t} \mu(|f|),\ \ \mu{\text -a.e.},\ t\ge 0,\end{align*} where  $c_0:=\|P_{t_1}\|_{L^1(\mu)\to L^\infty(\mu)}.$
This implies the desired assertion for $t_0>0$ such that $c_0^2\e^{-\ll t_0} \le \ff 1 2$.

(c)  Finally, combining    \eqref{AS}, \eqref{H}, \eqref{EP}    and \eqref{P}, we obtain
\beg{align*} &\|\nn_\GG P_{t+2 t_0 }f\|_{L^\infty(\mu)}^2  \le c_1 \| P_{t_0} |\nn_\GG  P_{t+t_0}f|^2 \|_{L^\infty(\mu)}\le 2c_1 \EE(P_{t+t_0}f, P_{t+t_0}f)\\
&\le c_2 \e^{-\ll t} \EE(P_{t_0}f, P_{t_0}f)
 \le c_2\eta_q(t_0) \e^{-\ll t} \mu\big((P_{t_0} |\nn_\GG f|^q)^{\ff 2 q}\big)\\
   &\le  c_3\eta_q(t_0) \e^{-\ll t} {\rm ess}_\mu\inf(P_{t+2t_0} |\nn_\GG f|^q)^{\ff 2 q} \end{align*}
for some constants $c_1, c_2, c_3 >0$. Then \eqref{GE} holds for $t_2= 2 t_0.$
\end{proof}

\section{Proof of Theorem \ref{T1.1} }
The first assertion is a generalization of the main result in \cite{LW} where $M=\R^d$ is considered.
As in \cite{LW}, the key point of the proof is to construct a coupling by parallel transform for long distance but by reflection for short distance.
The only difference is that we are working on a non-flat Riemannian manifold for which the curvature term appears in calculations.
Since It\^o's formula of the distance process has been well developed for couplings by both parallel displacement and reflection, the proof is also straightforward.

The proofs of the other two assertions are based on the log-Sobolev inequality and the log-Harnack inequality  derived in \cite{RW} and \cite{W10} respectively
for bounded below $\Ric_Z$.

\beg{proof}[Proof of Theorem \ref{T1.1}] (a) For two different points $x,y\in M$, let $P_{x,y}: T_xM \to T_yM$ be the parallel displacement along the minimal geodesic $\gg: [0,\rr(x,y)]\to M$ from $x$ to $y$, and let
$M_{x,y}:= P_{x,y}  - 2 \<\cdot, \dot \gg_{0}\> \dot\gg_{\rr(x,y)}: T_xM\to T_y M$ be the mirror reflection. Both  maps  are  smooth  in $(x,y)$ outside the cut-locus $\Cut(M)$. According to  \cite{K} and \cite{W94}, the appearance of the cut-locus and/or a convex boundary  helps for the success of coupling, i.e. it makes  the distance between two marginal processes smaller. So, for simplicity, we may and do assume that both the cut-locus and the boundary   are   empty,  see \cite[Section 3]{ATW06} or \cite[Chapter 2]{Wbook} for details.

Now,  let $X_t$ solve the SDE
$$\d_I X_t = \ss 2 u_t  \d B_t +Z(X_t)\d t,\ \ X_0=x,$$ where $\d_I$ denotes the It\^o differential introduced in \cite{E} on Riemannian manifolds,
$B_t$ is the $d$-dimensional Brownian motion, and $u_t$ is the horizontal lift of $X_t$ to the frame bundle $O(M)$.
Then $X_t$ is a diffusion process generated by $L$.  To construct the
 coupling by   reflection for short distance and parallel displacement  for long  distance, we introduce a
cut-off function $h\in C^1([0,\infty))$ which is decreasing such that  $h(r)=1$ for $r\le r_0,$   $h(r)=0$ for $r\ge r_0+1$, and $\ss{1-h^2}$ is also
in $C^1$, see e.g. \cite[(3.1)]{W15} for a concrete example. To construct the coupling in the above spirit, we split the noise into two parts, i.e. to replace $\d B_t$ by $h(\rr(X_t,Y_t))\d B_t' + \ss{1-h(\rr(X_t,Y_t))^2}\d B_t''$ for two independent Brownian motions $B_t'$ and $B_t''$,
 then make reflection for the $B_t'$ part and parallel displacement for the $B_t''$ part. More precisely, let $(X_t, Y_t)$ solve the following SDE  on $M\times M$ for $(X_0,Y_0)=(x,y)$:
\beg{align*} &\d_I X_t = \ss 2  \Big(h(\rr(X_t,Y_t)) u_t \d B_t' + \ss{1-h(\rr(X_t,Y_t))^2}u_t \d B_t''\Big) +Z(X_t)\d t,\\
&\d_I Y_t= \ss 2 \Big(h(\rr(X_t,Y_t))M_{X_t,Y_t} u_t \d B_t' + \ss{1-h(\rr(X_t,Y_t))^2}P_{X_t,Y_t}u_t \d B_t''\Big)+Z(Y_t)\d t.  \end{align*} Since the coefficients of the SDE are at least $C^1$ outside the diagonal $\{(z,z): z\in M\}$, it has a unique solution up to the coupling time
$$T:=\inf\{t\ge 0: X_t=Y_t\}.$$ We then let $X_t=Y_t$ for $t\ge T$ as usual. By the second variational formula and the index lemma (see e.g. the proof of \cite[Lemma 2.3]{WAnn}  and \cite[(2.4)]{W94}), the process $\rr_t:=\rr(X_t,Y_t)$ satisfies
$$\d \rr_t\le 2\ss 2 h(\rr_t)\d b_t + I_Z(X_t,Y_t)\d t,\ \ t\le T$$ for some one-dimensional Brownian motion $b_t$.  Thus, by condition \eqref{EB'},
\beq\label{ITO} \d \rr_t \le 2\ss 2 h(\rr_t)\d b_t +\big\{(K_1+K_2)1_{\{\rr_t\le r_0\}}-K_2\big\}\rr_t\d t,\ \ t\le T.\end{equation}
Since $h(\rr_t)=0$ for $\rr_t\ge r_0+1$ while  $\d\rr_t<0$ when $\rr_t\ge r_0+1,$  this implies
\beq\label{DD0}\rr_t\le (r_0+1)\lor\rr_0\le 1+r_0+\rr(x,y).\end{equation} On the other hand,  since $h(\rr_t)=1$ for $\rr_t\le r_0$, as observed in \cite{LW}
we have
\beq\label{DD}\E \rr_t\le c\e^{-\ll t}\rr(x,y),\ \ t\ge 0\end{equation} for some constants $c,\ll>0$. Indeed, let
$$\bar \rr_t= \vv \rr_t + 1- \e^{-N\rr_t},\ \ N=\ff {r_0} 2(  K_1+K_2), \vv=N\e^{-Nr_0}.$$
Then
$$\vv \rr_t\le \bar\rr_t\le (N+\vv)\rr_t,\ \ \ff{4N^2}{r(\vv\e^{Nr}+N)}\ge K_1+K_2 \ \text{for}\ r\in (0,r_0],$$ so that
 \eqref{ITO} and It\^o's formula lead to
\beg{align*}\d \bar\rr_t &\le 2\ss 2 (\vv+N\e^{-N\rr_t})h(\rr_t)\d b_t \\
&\qquad + (\vv+N\e^{-N\rr_t})\Big\{ (K_1+K_2)1_{\{\rr_t\le r_0\}}- K_2
-\ff {4N^2} {\rr_t (\vv\e^{N\rr_t}+N}1_{\{\rr_t\le r_0\}}\Big\}\rr_t\d t\\
&\le 2\ss 2 (\vv+N\e^{-N\rr_t})h(\rr_t)\d b_t -c_1\bar\rr_t\d t,\ \ t\le T\end{align*}
for some constant  $c_1$. This implies $\E \bar\rr_t\le  \bar\rr_0\e^{-c_1 t}.$ Then \eqref{DD} holds for some constants $c,\ll>0$.
Combining \eqref{DD0} with \eqref{DD} we arrive at
$$\E   \Phi(\rr_t/r) \le \sup_{s\in (0,  1+r_0+\rr_0]} \ff{\Phi(s/r)}{s} \E\rr_t \le c \e^{-\ll t} \rr(x,y) \sup_{s\in (0, 1+r_0+\rr_0]} \ff{\Phi(s/r)}{s}.$$ So,
\beg{align*}&W_{\Phi}(\dd_xP_t,  \dd_yP_t) \le \|\rr_t\|_{L^\Phi(\P)}= \inf\big\{r>0: \E \Phi(\rr_t/r)\le 1\big\}\\
&\le \inf\Big\{r>0:  \sup_{s\in (0, 1+ \rr(x,y)]} \ff{\Phi(\ff s r)}{s}\le \ff   {\e^{\ll t}}{c\rr(x,y)}\Big\},\end{align*}
which proves  \eqref{EB1'}. Therefore,  the proof of (1) is finished since   the second inequality therein is a simple consequence of \eqref{EB1'}.

(b) According to the proofs of \cite[Proposition 3.1 and Theorem 1.1]{WAnn}, our conditions imply that $P_t$ is hyperbounded; that is,  $\|P_t\|_{2\to 4}<\infty$ holds for some $t>0$. Since
\eqref{EB'} implies $\Ric_Z\ge -(K_1+K_2)$,
by the hyperboundedness and   \cite[Theorem 2.1]{RW},  we have  the defective log-Sobolev inequality
$$\mu(f^2\log f^2)\le C_1 \mu(|\nn f|^2) + C_2,\ \ f\in C_b^1(M), \mu(f^2)=1$$  for some constants $C_1, C_2>0$. Since  the symmetric Dirichlet form
$\EE(f,g):= \mu(\<\nn f,\nn g\>)$ with domain $H^{1,2}(\mu)$ is irreducible, according to \cite{W14} (see  also \cite{M}), the  log-Sobolev inequality \eqref{LS}   holds for some constant $C>0$, so that (2) is proved.

(c) According to \cite[Theorem 1.10]{S} (see   \cite{BGL0,W04,OV} for the case without boundary), the log-Sobolev inequality implies the Talagrand inequality
\beq\label{TL}W_2(f\mu, \mu)^2\le \ff C 2 \mu(f\log f),\ \ f\ge 0, \mu(f)=1. \end{equation} Next, let $P_t^*$ be the adjoint of $P_t$ in $L^2(\mu)$.
By Proposition \ref{PN}  for $P_t^*$ in place of $P_t$,  the log-Sobolev inequality   implies
\beq\label{EX0}\mu((P_t^*f)\log P_t^*f)\le \e^{-4t/C} \mu(f\log f),\ \ t\ge 0, f\ge 0, \mu(f)=1.\end{equation}
Moreover,   according to \cite[Theorem 1.1]{W10},   the curvature condition $\Ric_Z\ge -( K_1+K_2)=:-K$ is equivalent to the log-Harnack inequality
$$P_t (\log f )(x) \le \log P_tf (y) + \ff{K\rr(x, y)^2}{2(1-\e^{-2Kt})},\ \ t\ge 0, x,y\in M, 0\le f\in \B_b(M).$$
 By \cite[Proposition 1.4.4(3)]{WBook}, this implies
\beq\label{ET}\mu((P_t^* f)\log P_t^* f)\le \ff{K}{2(1-\e^{-2Kt})}W_2(f\mu,\mu)^2,\ \ f\ge 0, \mu(f)=1, t>0.\end{equation}
Combining \eqref{TL}, \eqref{EX0} and \eqref{ET}, we obtain
\beq\label{SL}\beg{split}& W_2((f\mu)P_{1+t}, \mu)^2= W_2((P_{1+t}^* f)\mu,\mu)^2\le \ff C 2 \mu((P_{1+t}^*f)\log P_{1+t}^*f)\\
&\le \ff C 2 \e^{-4t/C} \mu((P_1^* f)\log P_1^*f)
 \le c_1\e^{-4t/C} W_2(f\mu, \mu)^2,\ \ t\ge 0, f\ge 0, \mu(f)=1 \end{split}\end{equation} for some constant $c_1>0$. Noting that $\Ric_Z\ge -K$ implies $|\nn P_t f|\le \e^{Kt}P_t|\nn f|$
(see e.g. \cite{W10}), by Proposition \ref{T3.1} we have
$$W_2((f\mu)P_t,\mu)=W_2((f\mu) P_t, \mu P_t)\le c_2 W_2(f\mu,\mu),\ \ t\in [0,1], f\ge 0, \mu(f)=1.$$ Combining with \eqref{SL} yields
$$W_2( (f\mu) P_t,\mu)\le  c\e^{-\ll t} W_2(f\mu,\mu),\ \ t\ge 0, f\ge 0, \mu(f)=1$$ for some constants $c,\ll>0$. Therefore, the proof of  (3) is finished.
\end{proof}

\section{Proof of Theorem \ref{T1.2} and Corollary \ref{C1.3}}

\beg{proof}[Proof of Theorem \ref{T1.2}] (1) Since $\Ric_Z\ge -K$ for some constant $K\ge 0$, we have (see e.g. \cite{W10})
$$|\nn P_t f|\le \e^{Kt}P_t|\nn f|,\ \ f\in C_b^1(M).$$ Combining this with  Proposition \ref{P2.1} for $q=1$ and noting that $P_t|\nn f|$ is continuous, we obtain
$$|\nn P_t f| \le c_0\e^{-\ll t}  P_t|\nn f|,\ \ t\ge t_0, f\in C_b^1(M)$$ for some constants $c_0,\ll,t_0>0$. Obviously, \eqref{*NP} implies
$$\|\cdot\|_{L^1(P_t)}\le  \ff{\|\cdot\|_{L_*^\Phi(P_t)}}{\Phi^{-1}(1)},\ \ \Phi\in \bar {\scr N}.$$   Then
$$|\nn P_t f| \le \ff{c_0}{\Phi^{-1}(1)}\e^{-\ll t} \|\nn f\|_{L_*^\Phi(P_t)},\ \ t\ge 0, \Phi\in \bar{\scr N}, f\in C_b^1(M).$$
  According to
Proposition \ref{T3.1}, this is equivalent to
\beq\label{WW1} W_{\Phi}( \dd_x P_t,  \dd_yP_t)\le c_0 \Phi^{-1}(1) \e^{-\ll t} \rr(x,y),\ \ t\ge 0, x,y\in M. \end{equation} On the other hand,
noting that
$$\scr C( \dd_xP_t,  \dd_yP_t)\ni \pi_t:=( \dd_x P_t)\times ( \dd_yP_t) \le \|P_t\|_{L^1(\mu)\to L^\infty(\mu)}^2 (\mu\times \mu),$$
we obtain
$$W_{\Phi}( \dd_x P_t,  \dd_y P_t)\le \|\rr\|_{L^{\Phi}(\pi_t)} \le G_\Phi(2t),\ \ t>0.  $$
  Combining this with \eqref{WW1} and the semigroup property,  we arrive at
\beg{align*} W_{\Phi}( \dd_xP_t,  \dd_yP_t) \le\ff{c_0}{\Phi^{-1}(1)}\e^{-\ll t/2} W_\Phi(\dd_x P_{t/2}, \dd_y P_{t/2})
 \le \ff{c_0}{\Phi^{-1}(1)}\e^{-\ll t/2} G_\Phi(t).\end{align*}  This together with \eqref{WW1} implies \eqref{LL0}   for some constants $c,\ll>0.$ Moreover, \eqref{LL1} follows from \eqref{LL0} according to Proposition \ref{T3.1}.

 (2) By Proposition \ref{T3.1}, \eqref{WIF} implies
 $$|\nn P_t f|\le c \e^{-\ll t} P_t |\nn f|,\ \ \ t\ge 0, f\in C_b^1(M).$$ Then using the standard semigroup calculation of Bakry-Emery, this implies
 \beg{align*} &P_t (f^2\log f^2) - (P_tf^2)\log P_t f^2 = \int_0^t \ff{\d}{\d s} P_s\big\{(P_{t-s} f^2)\log P_{t-s} f^2\big\}\d s \\
 &= \int_0^t P_s\Big(\ff{|\nn P_{t-s} f^2|^2}{P_{t-s} f^2}\Big) \d s \le 4 c^2 \int_0^t \e^{-2\ll(t-s)} P_s \Big(\ff{(P_{t-s} \{f|\nn f|\})^2}{P_{t-s} f^2}\Big)\d s \\
 &\le 4 c^2 \int_0^t \e^{2\ll (t-s)} (P_t |\nn f|^2)\d s = \ff{2c^2 (1-\e^{-2\ll t})}{\ll} P_t |\nn f|^2,\ \ t\ge 0.\end{align*}
 Since $\lim_{t\to\infty} P_t g= \mu(g)$ for $g\in \B_b(M)$  due to the ergodicity, by letting $t\to\infty$ we prove the log-Sobolev inequality for
 \eqref{LS} for $C= \ff{2c^2}\ll.$
\end{proof}

\beg{proof}[Proof of Corollary \ref{C1.3}] We first observe that  the proof of  \cite[Theorem 4.2]{WAnn} works also for the non-symmetric case
with   $\nn Z$ in place of $\Hess_V$,  so that
\beq\label{ULT} \|P_t\|_{L^2(\mu)\to L^\infty(\mu)} \le \exp\Big[c + \ff c t\Big(1+\LL_1^{-1}(ct^{-1})+ \LL_2^{-1}(c^{-1}t)\Big)\Big],\ \ t>0.\end{equation}
Since in the symmetric case we have $\|P_t\|_{L^1(\mu)\to L^\infty(\mu)}\le \|P_{t/2}\|_{L^2(\mu)\to L^\infty(\mu)}^2$, the first assertion  follows immediately from Theorem \ref{T1.2}.

As for the non-symmetric case, since $$\|P_t\|_{L^1(\mu)\to L^\infty(\mu)}\le \|P_{t/2}\|_{L^1(\mu)\to L^2(\mu)}\|P_{t/2}\|_{L^2(\mu)\to L^\infty(\mu)},$$ by Theorem \ref{T1.2} and \eqref{ULT}
it suffices to prove
\beq\label{UU}  \|P_t\|_{L^1(\mu)\to L^2(\mu)}\le  c'+c'H(c't^{-1}),\ \ t>0\end{equation} for some constant $c'>0.$ According to \cite[Theorem 2.1]{RW},
\eqref{ULT} implies the super log-Sobolev inequality \eqref{SLS}   for
$$\bb(r):= c+ \ff c r\Big\{1+\LL_1^{-1}(cr^{-1})+ \LL_2^{-1}(c^{-1}r)\Big\},\ \ r>0$$ for some (possibly different) constant $c>0$. Then Proposition \ref{PN} with $p=1,q=2$ and
$\gg(r):=\ff{trh(r-1)}{(r-1)\int_0^1 s^{-1}h(s) \d s}$ implies
\eqref{UU}.

\end{proof}

\section{Proofs of Theorems  \ref{ST1}-\ref{ST2} and Proposition \ref{PPN}}

\beg{proof}[Proof of Theorems \ref{ST1}] Let $X_t(x)$ solve \eqref{SDE} with initial point $x$. By It\^o's formula and condition \eqref{DSS} we obtain
\beg{align*}&\d |X_t(x)-X_t(y)|^p\\
&\le \d M_t + p   |X_t(x)-X_t(y)|^{p-2}\bigg\{\ff{(p-2)|(\si(X_t(x))-\si(X_t(y))^*(X_t(x)-X_t(y))|^2}{|X_t(x)-X_t(y)|^2}\\
&\quad  +\|\si(X_t(x))-\si(X_t(y))\|^2_{HS}
+ 2 \<b(X_t(x)-b(X_t(y)), X_t(x)-X_t(y)\>\bigg\}\d t\\
&\le \d M_t - p K_p|X_t(x)-X_t(y)|^p \d t\end{align*}
for some martingale $M_t$. This implies
$$\E|X_t(x)-X_t(y)|^p\le \e^{-pK_p t}|x-y|^p,\ \ t\ge 0, x,y\in \R^d,$$ and thus,
\beq\label{LGR} \beg{split} |\nn P_tf(x)| &\le \limsup_{y\to x} \E\Big(\ff{|f(X_t(x))-f(X_t(y))|}{|X_t(x)-X_t(y)|}\cdot \ff{|X_t(x)-X_t(y)|}{|x-y|}\Big)\\
&\le \e^{-K_pt} (P_t |\nn f|^{\ff p{p-1}})^{\ff{p-1}p}.\end{split}\end{equation}
Then the desired assertion follows from Proposition \ref{T3.1}.
\end{proof}

\beg{proof}[Proof of Theorem \ref{ST2}] (1) We reformulate \eqref{SDE} as
\beq\label{SDE'} \d X_t= b(X_t)\d t +\ss 2\big(\si_0(X_t)\d B_t' + \ll_0\d B_t''\big),\end{equation} where $B_t'$ and $B_t''$ are independent  $d$-dimensional Brownian motions.
For any $x\ne y$, let $X_t$ solve this SDE with $X_0=x$, and let $Y_t$ solve the following coupled SDE with $Y_0=y$:
$$ \d Y_t= b(Y_t)\d t +\ss 2 \,\si_0(Y_t)\d B_t' + \ll_0\ss 2\,\bigg(\d B_t''- 2 \ff{\<X_t-Y_t, \d B_t''\>(X_t-Y_t)}{|X_t-Y_t|^2}\bigg).$$ That is, under the flat metric we have made coupling by reflection for $B_t''$ and coupling by parallel displacement for $B_t'$.
Obviously, the coupled SDE has a unique solution up to the coupling time
$$T_{x,y}:= \inf\{t\ge 0: X_t=Y_t\}.$$ We set $Y_t=X_t$ for $t\ge T_{x,y}$ as usual. Then by \eqref{DSS2} and It\^o's formula, we obtain
\beq\label{MM}  \d |X_t-Y_t|\le \d M_t + \big\{(K_1+K_2)1_{\{|X_t-Y_t|\le r_0\}} - K_2\big\} |X_t-Y_t|\d t,\ \ t\le T_{x,y}\end{equation}  for
$$\d M_t:=  \ff{\ss 2\<2\ll_0  \d B_t'' + (\si_0(X_t)-\si_0(Y_t))\d B_t', X_t-Y_t\> }{|X_t-Y_t|}$$
being  a martingale with
\beq\label{MM'}\d \<M\>_t\ge 8\ll_0^2\d t.\end{equation}  By repeating the argument leading to \eqref{DD}, it is  easy see that \eqref{MM} and \eqref{MM'}  imply
$$\E |X_t-Y_t|\le c\e^{-\ll t} |x-y|,\ \ t\ge 0$$ for some constants $c,\ll>0$ independent of $x,y$. Therefore,
$$|\nn P_t f|\le c \e^{-\ll t} \|\nn f\|_\infty,\ \ t\ge 0, f\in C_b^1(\R^d),$$ so that the first assertion follows from Proposition \ref{T3.1}.

(2) According to \cite[Theorem 1.1]{W11}, $a\ge \aa I$ and \eqref{DSS2'} imply the log-Harnack inequality
$$P_t (\log f )(x) \le \log P_tf (y) + \ff{c_1|x-y|^2}{1-\e^{-c_2t}},\ \ t\ge 0, x,y\in \R^d, 0\le f\in \B_b(\R^d)$$   for some constants $c_1,c_2>0$.
Combining this with the log-Sobolev inequality, we prove the second assertion as in (c) in the proof of  Theorem \ref{T1.1}.

(3) According to the proof of Theorem \ref{ST1}, the condition \eqref{DSS} implies the gradient estimate \eqref{LGR}. Next, by Proposition \ref{P2.1}, the ultracontractivity and \eqref{LGR} imply
 $$|\nn P_t f|\le c(p) \e^{-\ll t} (P_t |\nn f|^{\ff p{p-1}})^{\ff{p-1}p},\ \ t\ge 0, f\in C_b^1(\R^d)$$ for some $c(p)>0$ and $\ll>0$ independent of $p$. Then the proof if finished by Proposition \ref{T3.1}.
\end{proof}

\beg{proof}[Proof of Proposition \ref{PPN}] We will apply   results in  \cite{RW} and  \cite{WAD}. To this end, we introduce the Riemannian metric
$$g(\pp_i,\pp_j)= g_{ij}:= (a^{-1})_{i,j},\ \ 1\ge i,j\le d,$$  and let $\DD^g, \nn^g, \Hess^g$ be the corresponding Laplacian, gradient and Hessian tensor respectively. Then $L=\DD^g+Z$ for some $C^1$ vector field $Z$.   We first verify the Bakry-Emery curvature condition
\eqref{C} for some constant $K$. Using the Christoffel symbols, the intrinsic Hessian tensor induced by $g$ is formulated as
$$\Hess_f^g(\pp_i,\pp_j)= \pp_{ij}^2 f - \sum_{k=1}^d \GG_{ij}^k\pp_k f.$$
So, for any $x\in \R^d$ and $f\in C^2(\R^d)$ with $\Hess^g_f(x)=0$,   we have
$$\pp_{ij}^2f(x)= \sum_{n=1}^d \GG_{ij}^n \pp_n f(x),\ \ 1\le i,j\le d.$$
 Thus, by   Bochner-Weitzenb\"ock formula and \eqref{CC0},   at point $x$ there holds
 \beg{align*} &\Ric_Z(\nn^g f,\nn^g f) + K_0|\nn f|^2 = \ff 1 2 L\<a\nn f, \nn f\>- \<a\nn f, \nn Lf\>+K_0|\nn f|^2\\
 &\ge \ff 1 2\sum_{i,j,k,l=1}^d a_{kl} \Big[(\pp^2_{kl} a_{ij}) (\pp_if)(\pp_jf)+ 2 a_{ij} (\pp_{ki}^2 f) (\pp_{il}^2f)
 +2 (\pp_l a_{ij})\big\{(\pp^2_{ki}f)\pp_j f- (\pp_{ij}^2f)\pp_k f\big\}\Big]  \\
& = \ff 1 2 \sum_{i,j,k,l=1}^da_{kl} \Big[ (\pp^2_{kl} a_{ij}) (\pp_if)(\pp_jf)  +2 (\pp_l a_{ij})\sum_{n=1}^d (\pp_n f)\big\{\GG_{ki}^n \pp_i f- \GG_{ij}^n \pp_k f\big\}\Big] \\
&\ge - K_1|\nn f|^2 \ge -\ff{K_1}\aa \<a\nn f,\nn f\>=-\ff{K_1}\aa g(\nn^g f,\nn^gf) \end{align*} for some constant $K_1$.  Then \eqref{C} hold for some constant $K$.

Next,   \eqref{UL} implies that $P_t$ has a unique invariant probability measure $\mu$
 such that $\mu(\e^{c_2 |\cdot|^2})<\infty$ for some $c_2>\ff K {2\aa}$. By our assumption on $a$, the Riemannian distance $\rr $ induced by the metric $g$   is equivalent to the Euclidian metric:
\beq\label{EQ} \ff 1 {\|a\|_\infty} |\cdot|^2\le \rr_a^2(0,\cdot)\le\ff 1 \aa |\cdot|^2.\end{equation}  Then we may repeat
  the proof of \cite[Corollary 2.5]{RW} with $\gg(r)= c_2 r^{\dd} $ and $\rr=|\cdot|$ to prove
\beq\label{*D5} \|P_t\|_{L^2(\mu)\to L^\infty(\mu)} \le \exp\Big[c_3 t^{-\ff{\dd}{\dd-1}}\Big],\ \ t>0\end{equation} for some constant $c_3>0.$
Combining this with the curvature condition \eqref{C},  we obtain from \cite[Theorem 2.1]{RW} for $p=2$ and $q=\infty$ that
$$\mu(f^2\log f^2)\le r \EE(f,f) + c_4 r^{-\ff{\dd}{\dd-1}},\ \ r\in (0,1), \mu(f^2)=1$$
holds for some constant $c_4>0$. Applying Proposition \ref{PN} below for $p=1, q=2$ and $\gg(r)= c_5 t (r-1)^{\ff{\dd-1}{2\dd}-1}$ for constant $c_5>0$ such that $t= \int_1^2 \ff{\gg(r)}r\d r$, we obtain
$$\|P_t\|_{L^1(\mu)\to L^2(\mu)}\le \exp\Big[c_6 t^{-\ff{\dd}{\dd-1}}\Big],\ \ t\in (0,1)$$ for some constant $c_6>0$. Combining this with \eqref{*D5} we arrive at
$$\|P_t\|_{L^1(\mu)\to L^\infty(\mu)}\le c_7\exp\Big[c_7 t^{-\ff{\dd}{\dd-1}}\Big],\ \ t>0$$ for some constant $c_7>0$.

\end{proof}

\paragraph{Acknowledgement.} The author would like to thank  Jian Wang   for helpful comments.


\begin{thebibliography}{999}

\bibitem{AMS} S. Aida, T. Masuda, I. Shigekawa, \emph{Logarithmic
Sobolev inequalities and exponential integrability,} J. Funct.
Anal. 126(1994), 83--101.

\bibitem{ATW06} M. Arnaudon, A. Thalmaier, F.-Y. Wang, \emph{Harnack inequality and heat kernel estimate on manifolds with curvature unbounded below,} Bull. Sci. Math. 130(2006), 223--233.

 \bibitem{BE} D. Bakry,   M. Emery, \emph{Hypercontractivit\'e de
semi-groupes de diffusion}, C. R. Acad. Sci. Paris. S\'er. I Math.
299(1984), 775--778.
\bibitem{BGL0}  D. Bakry, I. Gentil, M. Ledoux, \emph{ Hypercontractivity of Hamilton-Jacobi equations,}
J. Math. Pures Appl.  80(2001), 669--696.
\bibitem{BGL} D. Bakry, I. Gentil, M. Ledoux, \emph{Analysis and Geometry of Markov Diffusion Operators,} Springer, 2014.

\bibitem{BGL15} D. Bakry, I. Gentil, M. Ledoux, \emph{On Harnack inqualities and optimal transportation,} Ann. Sc. Norm. Super. Pisa Cl. Sci. (5) XIV(2015), 705--727.
\bibitem{CW} M.-F. Chen, F.-Y. Wang, \emph{Estimates of logarithmic Sobolev constant: an improvement of Bakry-Emery criterion,} J. Funct. Anal. 144(1997), 287--300.
\bibitem{CW2} M.-F. Chen, F.-Y. Wang, \emph{Estimation of spectral gap for elliptic operators,} Trans. Amer. Math. Soc. 349(1997),
1239--1267.

\bibitem{DS} E. B. Davies, B. Simons, \emph{Ultracontractivity and heat kernel for Schr\"ordinger operators and Dirichlet Laplacians,} J. Funct. Anal.     59(1984), 335--395.


\bibitem{EB} A. Eberle, \emph{Reflection coupling and Warsserstein contractivity without convecity,} C. R. Acad. Sci. Paris. S\'er. I Math.
349(2011), 1101--1104.

\bibitem{EB2} A. Eberle, \emph{Reflection couplings and contraction rates for diffusions,} to appear in  Probab. Thero. Relat. Fields.

\bibitem{E} M. Emery, \emph{Stochastic Calculus in Manifolds,} Springer-Verlag, Berlin, 1989, with an appendix by P.-A. Meyer.

\bibitem{Gross} L. Gross, \emph{Logarithmic Sobolev inequalities
and contractivity properties of semigroups}, Lecture Notes in
Math. 1563, Springer-Verlag, 1993.

\bibitem{K}  W. S. Kendall, \emph{Nonnegative Ricci curvature and the
Brownian coupling property,} Stochastics 19(1986), 111--129.

\bibitem{KK} K.  Kuwada, \emph{Duality on gradient estimates and Wasserstein controls,} J. Funct. Anal. 258 (2010),   3758--3774.

\bibitem{KK2} K. Kuwada, \emph{Gradient estimate for Markov kernels, Wasserstein control and Hopf-Lax formula,} In $``$Potential Theory and Its Related Fields", 61--80, RIMS K\^oky\^uroku Bessatsu, B43, Res. Inst. Math. Sci. (RIMS), Kyoto, 2013.

 \bibitem{LW} D. Luo, J. Wang, \emph{Exponential convergence in $L^p$-Warsserstein distance for diffusion processes without uniform dissipative drift,} to appear in  Math. Nachr.

 \bibitem{M} L. Miclo, \emph{On hyperboundedness and spectrum of Markov operators,}  Inven. Math. 200(2015), 311--343.

\bibitem{O} F. Otto, \emph{The geometry of dissipative evolution equations: the porous medium equation,} Comm. Part. Diff.
Equat. 26(2001), 101--174.

\bibitem{OV} F.  Otto, C. Villani, \emph{Generalization of an inequality
by Talagrand and links with the logarithmic Sobolev inequality,}
J. Funct. Anal.  173(2000), 361--400.

\bibitem{PW}  E. Priola,  F.-Y. Wang, \emph{Gradient estimates for diffusion semigroups with singular coefficients,}  J. Funct. Anal. 236(2006), 244--264.

\bibitem{RS} M.-K. von Renesse, K.-T. Sturm, \emph{Transport inequalities, gradient estimates, entropy, and Ricci curvature,} Comm. Pure Appl. Math. 58(2005), 923--940.

\bibitem{RW} M. R\"ockner, F.-Y. Wang, \emph{Supercontractivity and ultracontractivity for
(non-symmetric) diffusion semigroups on manifolds,} Forum Math. 15(2003), 893--921.

\bibitem{RW04} M. R\"ockner, F.-Y. Wang, \emph{Spectrum for a class of (nonsymmetric) diffusion operators,} Bull. London Math. Soc. 36(2004), 95--104.

\bibitem{S} J. Shao, \emph{Hamilton-Jacobi semi-groups in infinite dimensional
spaces,} Bull. Sci. math. 130(2006), 720--738.

\bibitem{ST}   W. Stannat,  \emph{(Nonsymmetric) Dirichlet operators on $L^1$: Existence, uniqueness and associated Markov processes,}
Ann. Sc. Norm. Super. Pisa Cl. Sci.   28(1999) 99--140.

\bibitem{V2} C. Villani, \emph{Topics in Optimal Transportation,} Amer. Math. Soc. 2003.

\bibitem{V} C. Villani, \emph{Optimal Transpot, Old and New,} Springer, 2009.

\bibitem{W94} F.-Y. Wang, \emph{Application of coupling method to the
Neumann eigenvalue problem,} Probab. Theory Relat. Fields
98(1994), 299--306.

\bibitem{W97} F.-Y. Wang, \emph{Logarithmic Sobolev inequalities on noncompact Riemannian manifolds,} Probab.
Theory Relat. Fields 109(1997), 417--424.
\bibitem{W00}  F.-Y. Wang, \emph{Functional inequalities for empty essential spectrum,}  J. Funct. Anal. 170(2000), 219--245.

\bibitem{W04} F.-Y. Wang, \emph{Probability distance inequalities on Riemannian manifolds and path spaces,} J. Funct. Anal. 206(2004), 167--190.
 \bibitem{Wbook}  F.-Y. Wang, \emph{Functional Inequalities, Markov Semigroups and Spectral Theory,} Science Press, 2005, Beijing.

 \bibitem{WAnn} F.-Y.Wang, \emph{Log-Sobolev inequalities: different roles of Ric and Hess,}  Ann. Probab. 37(2009), 1587--1604.

 \bibitem{WAD} F.-Y.Wang, \emph{Log-Sobolev inequality on non-convex manifolds,}   Adv. Math. 222(2009), 1503--1520.

 \bibitem{W10} F.-Y.Wang, \emph{Harnack inequalities on manifolds with boundary and applications,} J. Math. Pures Appl. 94(2010) 304--321.

\bibitem{W11} F.-Y.  Wang, \emph{ Harnack inequality for SDE with multiplicative noise and extension to Neumann semigroup on nonconvex manifolds,} Ann. Probab. 39(2011),  1449--1467.

\bibitem{W14} 	F.-Y.Wang, \emph{Criteria on spectral gap of Markov operators,} J. Funct. Anal. 266(2014), 2137--2152.

\bibitem{WBook} F.-Y. Wang, \emph{Analysis for Diffusion Processes on Riamnnian Manifolds,}   Springer, 2014.
\bibitem{W15} F.-Y. Wang,  \emph{Asymptotic couplings by reflection and applications for
nonlinear monotone SPDEs,} Nonlinear Anal. 117(2015), 169--188.





\end{thebibliography}
\end{document}